\documentclass[11pt]{amsart}

\usepackage{amsmath, amsthm}
\usepackage{amssymb}
\usepackage{amscd}
\usepackage{palatino}
\usepackage{latexsym}
\usepackage{epsfig}
\usepackage{graphics}
\usepackage{amsfonts}
\usepackage{psfrag}
\usepackage{graphicx}

\input xy
\xyoption{all}
\UseComputerModernTips

\numberwithin{equation}{section}
\numberwithin{figure}{section}

\newtheorem{theorem}{Theorem}[section]

\newtheorem{proposition}[theorem]{Proposition}

\theoremstyle{definition}

\theoremstyle{remark}
\newtheorem{remark}[theorem]{Remark}
\newcommand{\Z}{{\mathbb{Z}}}
\newcommand{\Q}{{\mathbb{Q}}}
\newcommand{\R}{{\mathbb{R}}}
\newcommand{\N}{{\mathbb{N}}}

\newcommand{\g}{\mathfrak{g}}

\newcommand{\into}{\hookrightarrow}
\newcommand{\onto}{\twoheadrightarrow}

\renewcommand{\mod}{/\!/}

\DeclareMathOperator{\Crit}{Crit}
\DeclareMathOperator{\ev}{ev}
\DeclareMathOperator{\Map}{Map}

\DeclareMathOperator{\Hol}{Hol} 
\DeclareMathOperator{\hol}{hol} 
\DeclareMathOperator{\Fred}{Fred}

\newcommand{\hs}{{\hspace{3mm}}}
\newcommand{\hsm}{{\hspace{1mm}}}

\begin{document}

%%%%%%%%%%%%%%%%%%%%
% Disclaimer
%%%%%%%%%%%%%%%%%%%%
%\begin{center}
%\framebox{
%{\Large\bf DRAFT (\today): DO NOT DISTRIBUTE.}}
%\end{center}

\title{Kirwan surjectivity in $K$-theory for Hamiltonian loop group quotients}

\author{Megumi Harada} 
\address{Department of Mathematics and
Statistics, McMaster University, 1280 Main Street West, Hamilton,
Ontario L8S4K1, Canada} 
\email{Megumi.Harada@math.mcmaster.ca}
\urladdr{http://www.math.mcmaster.ca/{}Megumi.Harada/}

\author{Paul Selick} 
\address{Department of Mathematics, University of Toronto, 
40 St. George Street, Toronto, Ontario M5S2E4, Canada}
\email{selick@math.toronto.edu}

\keywords{Hamiltonian actions, loop groups, symplectic quotients, $K$-theory}
\subjclass[2000]{Primary: 53D20; Secondary: 55N15}

\date{\today}

%%%%%%%%%%%%%%%%%%%%%
%  Abstract
%%%%%%%%%%%%%%%%%%%%%

\begin{abstract}

Let $G$ be a compact Lie group and $LG$ be its associated loop group. 
The main result of this manuscript is a surjectivity theorem
from the equivariant $K$-theory of a Hamiltonian $LG$-space 
onto the integral $K$-theory of its Hamiltonian $LG$-quotient. 
Our result is
a $K$-theoretic analogue of previous work in rational Borel-equivariant
cohomology by Bott, Tolman, and Weitsman. Our proof techniques differ from that of Bott, 
Tolman, and Weitsman in that they explicitly use the Borel construction, which we do not 
have at our disposal in equivariant $K$-theory; we instead directly construct $G$-equivariant 
homotopy equivalences to obtain the necessary isomorphisms in equivariant $K$-theory. 
The main theorem should also be viewed as a first step toward a similar theorem in $K$-theory for 
quasi-Hamiltonian $G$-spaces and their associated quasi-Hamiltonian quotients. 
\end{abstract}

\maketitle

\tableofcontents

\section{Introduction}\label{sec:intro}

The main result of this manuscript is a surjectivity theorem from the equivariant $K$-theory of a Hamiltonian loop group space onto the integral $K$-theory of its Hamiltonian quotient. Our result should be understood as an integral $K$-theoretic analogue of a result of Bott, Tolman, and Weitsman \cite[Theorem 2]{BTW04}, which is in turn a loop group analogue of Kirwan's surjectivity theorem in \cite{Kir84}. 

We now briefly recall the setting of Kirwan's original result, which is an influential tool in the theory of Hamiltonian spaces and their associated quotients. For $G$ a compact Lie group and $M$ a finite-dimensional Hamiltonian $G$-space with proper moment map \(\mu: M \to \g^{*},\) Kirwan shows that the norm-square of the moment map \(\|\mu\|^{2}: M \to \R\) is a $G$-equivariantly perfect Morse function on $M$, and hence concludes that the inclusion \(\mu^{-1}(0) \into M\) of the zero-level set of $\mu$ into $M$ induces a {\em surjection} in $G$-equivariant cohomology 
\begin{equation}\label{eq:Kirwan-original} 
H^{*}_{G}(M;\Q) \onto H^{*}_{G}(\mu^{-1}(0); \Q). 
\end{equation}
If the $G$-action on $\mu^{-1}(0)$ is free, the right hand side of~\eqref{eq:Kirwan-original} is isomorphic to the ordinary cohomology \(H^{*}(\mu^{-1}(0)/G; \Q)\) of the symplectic quotient \(M \mod G := \mu^{-1}(0)/G,\) and hence we obtain the Kirwan surjection \(\kappa: H^{*}_{G}(M;\Q) \onto H^{*}(M \mod G; \Q).\) The map $\kappa$ has been a powerful tool in equivariant symplectic geometry since its introduction in Kirwan's 1984 manuscript, since by explicitly computing \(\ker(\kappa)\) one may then obtain explicit descriptions of the cohomology rings of finite-dimensional symplectic quotients (see e.g. \cite{JK95, TW03, Gol01}). Moreover, the first author and Landweber have recently
generalized the surjectivity result~\eqref{eq:Kirwan-original} of Kirwan as well as related results for computing 
$\ker(\kappa)$ to the setting of topological $K$-theory \cite{HarLan07, HL-kernel}. Since integral $K$-theory encodes more delicate (e.g. torsion) information than does rational cohomology, these are promising new developments in the computation of the topology of Hamiltonian quotients. 

A natural analogue of the above finite-dimensional theory of Hamiltonian symplectic quotients exists when the compact Lie group $G$ is replaced by the infinite-dimensional Banach Lie group $L_{s}G$ consisting of loops \(\gamma: S^{1} \to G\) on $G$ of Sobolev class $s$ (details in Section~\ref{sec:background}). A Hamiltonian $L_{s}G$-space is then an infinite-dimensional Banach manifold ${\mathcal M}$ on which $L_{s}G$ acts with moment map \(\mu: {\mathcal M} \to L_{s}\g^{*}.\) We always assume $\mu$ is proper. The Hamiltonian $L_{s}G$-quotient is defined, when $G$ acts freely on $\mu^{-1}(0)$, to be ${\mathcal M} \mod L_{s}G := \mu^{-1}(0)/G$, just as in the finite-dimensional case. 
Many symplectic manifolds arise as Hamiltonian $L_{s}G$-quotients, the most famous class of which are the moduli spaces of flat $G$-connections over Riemann surfaces with boundary. Hence it is of importance to obtain explicit methods to compute topological invariants such as rational cohomology or integral $K$-theory corresponding to these quotients. 

As mentioned above, Bott, Tolman, and Weitsman already proved in \cite{BTW04} the loop group analogue of~\eqref{eq:Kirwan-original}, i.e. the Kirwan surjectivity in rational cohomology for loop group 
quotients.  Our main result is the following, which is in turn the integral $K$-theory version of Bott, Tolman, and Weitsman's result and therefore the first step toward explicit computations of the $K$-theory of Hamiltonian $L_{s}G$-quotients. 

\begin{theorem}\label{theorem:main} 
Let $G$ be a compact Lie group equipped with a bi-invariant metric. Let $L_{s}G$ be the corresponding loop group for \(s > \frac{3}{2}.\) Let \(({\mathcal M}, \omega, \mu)\) be a Hamiltonian $L_{s}G$-space with proper moment map \(\mu: {\mathcal M} \to L_{s}\g^{*}.\) Then the inclusion \(\mu^{-1}(0) \into {\mathcal M}\) induces a surjection in equivariant $K$-theory
\[
K^{*}_{G}({\mathcal M}) \onto K^{*}_{G}(\mu^{-1}(0)).
\]
In particular, if $G$ acts freely on $\mu^{-1}(0)$, the Kirwan map obtained by the composition 
\[
\kappa: K^{*}_{G}({\mathcal M}) \onto K^{*}_{G}(\mu^{-1}(0)) \stackrel{\cong}{\to} K^{*}({\mathcal M} \mod L_{s}G = \mu^{-1}(0)/G)
\]
is surjective. 
\end{theorem}

Before discussing the proof of our main theorem, we take a moment to discuss the generalized equivariant cohomology theory which is the focus of this manuscript. 
By $K$-theory, we mean topological, 
integral, complex $K$-theory, taking $K^0(X)$ to be the 
isomorphism classes of virtual complex vector bundles over $X$ when $X$ is compact, or 
for more general $X$, taking $[X,\Fred(\mathcal{H})]$ for a
complex separable Hilbert space $\mathcal{H}$. 
In the
equivariant case, by $K_G(X)$ we mean Atiyah-Segal $G$-equivariant $K$-theory
\cite{Seg68}, built from $G$-equivariant vector bundles if $X$ is a compact $G$-space, and 
$G$-equivariant maps 
$[X, \Fred(\mathcal{H}_G)]_G$ if $X$ is (as in our case) noncompact (here $\mathcal{H}_G$ contains
infinitely many copies of 
every irreducible representation of $G$, see e.g. \cite{AtiSeg04}). 

We now briefly describe our arguments in the proof of Theorem~\ref{theorem:main}. This also provides another motivation for our result, as it involves the theory of quasi-Hamiltonian $G$-spaces as introduced by Alekseev, Malkin, and Meinrenken \cite{AMM98}. A quasi-Hamiltonian $G$-space $M$ is a $G$-space equipped with a $2$-form $\omega$ and a quasi-Hamiltonian
moment map $\Phi$. The $2$-form $\omega$ is not necessarily symplectic
and the moment map $\Phi$ does not necessarily satisfy Hamilton's
equation, but both of these discrepancies from the classical Hamiltonian theory
are well-controlled by the canonical $3$-form on the 
Lie group $G$. We will not need much of the theory of quasi-Hamiltonian $G$-spaces in this manuscript so we will not present all details here; suffice it to mention that one of the essential observations (and indeed, the original motivation) for this theory is
that quasi-Hamiltonian $G$-spaces \((M, \omega, \Phi)\) are in one-to-one
correspondence with Hamiltonian $L_sG$-spaces with proper moment
map. There is a holonomy map  \(L_{s}\g^{*} \to G\) making $L_{s}\g^{*}$ a fiber bundle over $G$, and a theorem of Alekseev, Malkin, and Meinrenken \cite[Theorem 8.3]{AMM98} states that any Hamiltonian $L_{s}G$-space with proper moment map is obtained as a pullback of the fibration $L_{s}\g^{*} \to G$ via a quasi-Hamiltonian moment map \(\Phi: M \to G\) for some quasi-Hamiltonian space $M$. Moreover, the Hamiltonian $L_{s}G$-quotient ${\mathcal M} \mod L_{s}G$ may be identified with the quasi-Hamiltonian quotient \(M \mod G := \Phi^{-1}(e)/G.\) Thus the theory of Hamiltonian $L_{s}G$-quotients is inextricably linked with that of quasi-Hamiltonian quotients. Indeed, Bott, Tolman, and Weitsman prove in \cite[Theorem 3]{BTW04} a result in rational cohomology for quasi-Hamiltonian quotients similar in spirit to the original Kirwan surjectivity theorem, using their surjectivity theorem for Hamiltonian $L_{s}G$-quotients. We expect that a similar result should hold also in integral $K$-theory; hence this is another motivation for the current manuscript. 

The rational cohomology version \cite[Theorem 2]{BTW04} of our main theorem was proven by Bott, Tolman, and Weitsman by fully exploiting the above-mentioned relationship between an (infinite-dimensional) Hamiltonian $L_{s}G$-space ${\mathcal M}$ and its associated finite-dimensional quasi-Hamiltonian $G$-space $M$, and in particular the fact that ${\mathcal M}$ is a fiber product of $M$ with the space $L_{s}\g^{*}$. The essential technical point is the following: the main difficulties of dealing with Hamiltonian $L_{s}G$-spaces arise from their infinite-dimensionality, but it is precisely the Alekseev-Malkin-Meinrenken description of ${\mathcal M}$ as a fiber product of a finite-dimensional space $M$ with a specific infinite-dimensional affine space $L_{s}\g^{*}$ (about which a great deal is known) which makes these technical problems surmountable. Indeed, using the description of $L_{s}\g^{*}$ as an affine space of connections on a $G$-principal bundle over $S^{1}$, Bott, Tolman, and Weitsman use the classical infinite-dimensional Morse-theoretic approach of approximating continuous paths on $G$ by piecewise smooth geodesics to construct a sequence of finite-dimensional approximating spaces $Y_{n}$ to ${\mathcal M}$. Proving the surjectivity for each of these $Y_{n}$ then implies surjectivity for ${\mathcal M}$. In rational cohomology, this last step of proving the surjectivity for each $Y_{n}$ is accomplished in \cite{BTW04} by a careful analysis of the local geometry of a Morse-type function defined on each $Y_n$, in addition to the rational cohomology Atiyah-Bott lemma \cite[Prop 13.4]{AtiBot83}, \cite[Lemma 1.4]{BTW04}. In the setting of $K$-theory, the analysis of the local geometry of course remains the same but we must instead rely on the integral topological $K$-theory version of the Atiyah-Bott lemma, developed and used in \cite{HarLan07} (see also \cite{VeVi03} for a version in algebraic $K$-theory). 

Our arguments in the present manuscript follow in broad outline that given by Bott, Tolman, and Weitsman; in particular, we use the same approximating spaces as constructed in \cite{BTW04}. However, the key difference between our arguments and those in \cite{BTW04} is that Bott, Tolman, and Weitsman consistently use explicitly the Borel construction in equivariant cohomology, as well as other techniques specific to cohomology, none of which we have at our disposal in equivariant $K$-theory. Hence, by necessity, we must use other methods to prove that these approximating spaces do also indeed approximate ${\mathcal M}$ (in the appropriate sense) in $K$-theory.  We accomplish this by taking the arguments in \cite{BTW04} a few steps further and exhibiting explicit $G$-equivariant homotopy equivalences between the appropriate spaces, from which the relevant isomorphisms on the equivariant $K$-theories of course immediately follow. In other words, our proofs are based on exhibiting explicit geometric relationships on the underlying spaces, instead of algebraic arguments on the level of the associated algebraic invariants.

We close with some comments on directions for future work. Firstly, as was already mentioned above, the present manuscript should be viewed as a first step toward a generalization to $K$-theory of the Kirwan surjectivity-type theorem for quasi-Hamiltonian quotients proven by Bott, Tolman, and Weitsman in rational cohomology \cite[Theorem 3]{BTW04}. Secondly, it is evident that, given such a surjectivity theorem as in Theorem~\ref{theorem:main}, some explicit computations of the integral $K$-theory of specific examples (such as the moduli spaces of flat connections on $G$-bundles over Riemann surfaces) are clearly in order.  We intend to address these and related issues in future work.

The contents of this manuscript are as follows. In Section~\ref{sec:background}
we recall only the bare minimum of the necessary technicalities of our setup in order to be able to 
state our results. In Section~\ref{sec:infinite} we construct an
infinite-dimensional space $Y$ which approximates the original
Hamiltonian $L_{s}G$-space ${\mathcal M}$, in the sense that if we
prove an analogous surjectivity theorem for $Y$ then this implies our
main theorem. In Section~\ref{sec:finite} we construct a sequence
$Y_{n}$ of finite-dimensional spaces which in turn approximate $Y$, again in
the sense that if we prove the analogous surjectivity theorems
for each of the $Y_{n}$ (compatible in an appropriate sense) then we
may deduce the surjectivity for $Y$. This proves the main result.  We
have relegated to an Appendix some technical and straightforward
arguments in $G$-equivariant homotopy theory.

\medskip
\noindent{\bf Acknowledgements.}  We thank Jonathan Weitsman, Lisa Jeffrey, Eckhard Meinrenken, 
and Greg Landweber for their 
interest and encouragement. The first author thanks Andrew Nicas for many useful conversations.
The authors were supported in part by NSERC Discovery Grants.

\section{Background: Hamiltonian $L_sG$-spaces and quasi-Hamiltonian $G$-spaces}\label{sec:background}

We now turn to some of the technical analytical details of our
infinite-dimensional problem. We refer the reader to \cite{AMM98, BTW04} for
details.

Let $G$ be a compact Lie group equipped
with a bi-invariant metric.
By the loop group $L_sG$ we will mean
the infinite-dimensional space 
\[
L_sG := \left\{ \gamma: S^1 \to G \,\, \bigg\vert \,\, \gamma\hsm \hsm  \mbox{is of
  Sobolev class} \hsm \hsm s   \right\},
\]
for some fixed \(s > \frac{3}{2} .\) 
Then 
  $L_sG$ is itself a Banach Lie group 
where the group
  structure is given by pointwise multiplication in $G$. The Lie
  algebra $L_s\g$ of $L_sG$ may be identified with the
  space \(\Omega_s^0(S^1; \g)\) of maps from $S^1$ to $\g$ of
  Sobolev class $s$. We then define $L_s\g^*$ to be 
 \(\Omega^1_{s-1}(S^1; \g),\) the
  space of $\g$-valued one-forms on $S^1$ of Sobolev class
  $s-1$. The pairing 
\[
(\alpha, \beta) \in L_{s}\g \times L_{s}\g^{*} \mapsto \int_{S^{1}} \langle \alpha, \beta \rangle 
\]
where $\langle \cdot, \cdot \rangle$ is the given metric restricted to \(\g = T_{e}G\) induces an inclusion \(L_{s}\g^{*} \into (L_{s}\g)^{*}\) into the Banach space dual.  Viewing \(L_{s}\g^{*} := \Omega^{1}_{s-1}(S^{1}; \g)\) as the affine space of connections on the trivialized $G$-bundle \(S^{1} \times G \to S^{1},\) we let the loop group $L_{s}G$ act on $L_{s}\g^{*}$ via gauge transformations, namely for \(\gamma \in L_{s}G, \beta \in L_{s}\g^{*}\) we have 
\begin{equation}\label{eq:gauge-action}
\gamma \cdot \beta = \gamma \beta \gamma^{-1} - \gamma^{*} \overline{\theta},
\end{equation}
where $\overline{\theta}$ denotes (following notation of \cite{AMM98}) the right-invariant Maurer-Cartan form on $G$. We note in particular that this action of $L_sG$ on $L_s\g^*$ is not linear because of the presence of the translation term $\gamma^*\overline{\theta}$. Indeed, this action of $L_sG$ on $L_s\g^*$ is in fact induced from the action of an $S^1$-central extension $\widehat{L_sG}$ of $L_sG$ on its own Lie algebra dual $\widehat{L_s\g}^* \cong L_s\g^* \oplus \R$, where the second factor is the Lie algebra of the central $S^1$. Since the $S^1$ is central, the coadjoint action of $\widehat{L_sG}$ descends to an action of $L_sG$, and the codimension-$1$ affine subspaces ``at level $t$'' \(L_s\g^* \oplus \{t\} \subseteq L_s\g^* \oplus \R \cong \widehat{L_s\g}^*\) for \(t \in \R\) are preserved under this action. The gauge action~\eqref{eq:gauge-action} is then the restriction of this action to the ``level $1$'', i.e. we have identified $L_s\g^*$ with the codimension-$1$ affine subspace \(L_s\g^* \oplus \{1\} \subseteq L_s\g^* \oplus \R \cong \widehat{L_s\g}^*.\) (See e.g. \cite[Section 4]{PS86}.) Hence although by \(0 \in L_s\g^*\) we do mean the identically $0$ connection $1$-form in $\Omega^1_{s-1}(S^1; \g)$, we do not mean the zero element in the dual vector space $\widehat{L_s\g}^*$ (since this $0$ connection is actually identified with $(0,1) \in L_s\g^* \oplus \R$). In particular, the stabilizer group of this $0$ connection in $L_s\g^*$ is not all of $L_sG$ but only the subgroup of identically constant loops \(G \subseteq L_sG.\)

Continuing to interpret $L_{s}\g^{*}$ as a space of connections, we also have a holonomy map \(\hol: L_{s}\g^{*} \to P_{e}G^{*}\) where $P_{e}G^{*}$ denotes\footnote{For consistency we here follow the notation of \cite{BTW04}; however, we do point out that the $P_eG^*$ simply denotes a space of continuous loops, and is not a ``dual'' of any space.} the space of continuous paths in $G$ based at the identity \(e \in G\) (that $\hol(\beta)$ is indeed continuous follows from our choice of Sobolev parameter). 
We also let \(\Hol: L_s\g^* \to G\) denote the ``endpoint'' holonomy \(\Hol = \hol(1)\),  i.e. the holonomy around the loop \(S^{1} = \R/\Z.\) Recall also that \(\Hol: L_{s}\g^{*} \to G\) makes $L_{s}\g^{*}$ a principal $\Omega_{s}G$-bundle over $G$, where \(\Omega_{s}G\) is the subgroup of loops based at $e \in G$ (see \cite{AMM98}).  

A Hamiltonian $L_{s}G$-space is then defined to be a Banach manifold ${\mathcal M}$ equipped with an $L_{s}G$-action, an $L_{s}G$-invariant $2$-form $\omega$ and an equivariant map \(\mu: {\mathcal M} \to L_{s}\g^{*}\) such that the $2$-form $\omega$ is closed and is weakly non-degenerate (i.e. induces an injection \(T_{x}{\mathcal M} \to T^{*}_{x}{\mathcal M}\) for all \(x \in {\mathcal M}\)) and $\mu$ is a moment map for the $L_{s}G$-action. In \cite{AMM98}, Alekseev, Malkin, and Meinrenken prove that every Hamiltonian $L_{s}G$-space with a proper moment map \(\mu: {\mathcal M} \to L_{s}\g^{*}\) arises as a pullback $\Omega_{s}G$-fibration 
\begin{equation}\label{eq:pullback-diagram} 
\xymatrix{
{\mathcal M} \ar[r] \ar[d] & L_{s}\g^{*} \ar[d]^{\Hol} \\
M \ar[r]_{\Phi} & G 
}
\end{equation}
where $M$ is a finite-dimensional smooth manifold and $\Phi: M \to G$ a smooth map. In fact $M$ is naturally equipped with a $G$-action and a $2$-form, making the map $\Phi$ $G$-equivariant with respect to the given $G$-action on $M$ and the conjugation action of $G$ on itself.  Thus we may conclude that $M$ a quasi-Hamiltonian space in the sense of \cite[Definition 2.2]{AMM98}. We will not need further specific properties of quasi-Hamiltonian spaces here, so we do not pursue this matter in more detail. The essential fact we need is simply the existence of the above pullback diagram for some $G$-space $M$ and $G$-equivariant map $\Phi$.

\section{Surjectivity via an infinite-dimensional approximation}\label{sec:infinite}

In this section we construct an infinite-dimensional space which will
serve as an effective replacement, for the purposes of computing
equivariant $K$-theory of our
original Hamiltonian $L_sG$-space. The essential point is that we wish
to replace $L_s\g^*$, a space of connections on a principal
$G$-bundle, with a space $Y$ which is intimately related to a much
more familiar object in
classical infinite-dimensional Morse theory (and for which
Morse-theoretic tools have been well-developed), namely, the space of
piecewise smooth paths on $G$. Our construction of this space $Y$ will
be the same as that given in \cite{BTW04}, so we only briefly recall
it here. Our notation closely follows that in \cite{BTW04}. 

The main result of this section is the following. On this
infinite-dimensional approximating space $Y$ we will see that there exists a function
$\hat{f}: Y \to \R$ related to the classical energy functional on
a space of piecewise smooth paths. We prove below that the map induced
by the inclusion \(\mu^{-1}(0) \into {\mathcal M},\)
\begin{equation}\label{eq:Kirwan-to-mu-level}
K^*_G({\mathcal M}) \to K^*_G(\mu^{-1}(0))
\end{equation}
is surjective if and only if 
\begin{equation}\label{eq:Kirwan-to-hatf-level}
K^*_G(Y) \to K^*_G(\hat{f}^{-1}(0))
\end{equation}
is surjective, where the map~\eqref{eq:Kirwan-to-hatf-level} is also induced by the inclusion
\(\hat{f}^{-1}(0) \into Y.\) In other words, we reduce the main
problem, that of proving surjectivity of~\eqref{eq:Kirwan-to-mu-level}, to one of analyzing the infinite-dimensional Morse theory of
$Y$. We exploit this further in Section~\ref{sec:infinite}
by constructing a sequence of finite-dimensional approximations to $Y$
given by piecewise geodesic paths. 

Our proof of the equivalence of the surjectivity
of~\eqref{eq:Kirwan-to-mu-level} and~\eqref{eq:Kirwan-to-hatf-level}
does not follow that of \cite{BTW04}, since their arguments make
explicit use of the Borel construction and other techniques specific
to ordinary cohomology which in our setting we do not have at our disposal.  This
turns out not to be a serious limitation, since certain geometric
constructions given in \cite{BTW04} can, by some further argument, be
seen to be $G$-equivariant homotopy equivalences. This stronger
statement then implies the necessary results also in equivariant
$K$-theory.

We now briefly recall the construction of the infinite approximating space $Y$
to ${\mathcal M}$ as given in \cite[Section 3]{BTW04}. As already mentioned,
the
construction depends substantially on the description of a Hamiltonian
$L_sG$-space ${\mathcal M}$ as a pullback of an $\Omega_sG$-bundle over
$G$ as in \cite{AMM98}. 
 Let $M$ be the quasi-Hamiltonian
space associated to ${\mathcal M}$, so that \(\Omega_sG \into
{\mathcal M} \to M\) is the fibration in the left vertical
of~\eqref{eq:pullback-diagram}. Then 
\[
{\mathcal M} \cong \left\{ (\gamma, m) \in L_s\g^* \times M \, \vert
\, \Phi(m) = \Hol(\gamma) \right\}.
\]
The idea of the
construction of the space $Y$ is to replace the contractible space $L_s\g^*$
with a different contractible space which is related to $L_s\g^*$ by
taking holonomies. For instance, the holonomy map
$\hol$ on $L_s\g^*$ takes values in $P_eG^*$, the space of 
continuous maps (equipped with uniform topology) 
\(\gamma: [0,1] \to
G\) with \(\gamma(0) = e \in G,\) where $e$ is the identity of
$G$.  The group $G$ acts on
$P_eG^*$ by pointwise conjugation \((g\cdot \lambda)(t) := g^{-1}\lambda(t)g,\)
and $P_eG^*$ is $G$-equivariantly contractible. Let  \(\rho^*:
P_eG^* \to G\) denote the endpoint map \(\gamma \mapsto \gamma(1) \in G.\) We
may now define, analogously to ${\mathcal M}$, the space 
\[
{\mathcal M}^* := \left\{ (\gamma, m) \in  P_eG^* \times M \, \vert \,
\Phi(m) = \rho^*(\gamma) \right\}.
\]
The space ${\mathcal M}^*$ is also equipped with a $G$-action, given by the diagonal action. 
Projection to the
first factor yields a $G$-equivariant map \(\mu^*:
{\mathcal M}^* \to P_eG^*.\) If we denote by ${\mathbf e}$ the
constant path at the identity in $P_eG^*$, we note that
the holonomy
map \(\hol: L_s\g^* \to P_eG^*\) induces an inclusion \({\mathcal M}
\into {\mathcal M}^*\) which takes \(\mu^{-1}(0)\)
exactly to $(\mu^*)^{-1}({\mathbf e})$.

We have now related the question of analyzing the level set
\(\mu^{-1}(0) \subseteq {\mathcal M}\) to that of the level set 
$(\mu^*)^{-1}({\mathbf e})$ in the space ${\mathcal
  M}^*$. However, ${\mathcal M}^*$ is not the approximating space to
${\mathcal M}$ that we really want, most especially since (as we have
just seen) ${\mathcal M}^*$ is {\em larger} than ${\mathcal
  M}$. However, now that we have phrased a question in terms of
continuous paths in $G$, we may use the analytical trick of
approximating continuous functions by piecewise smooth ones. 

With this in mind, let \(\widehat{P_eG}\) denote the space of {\em
  piecewise smooth} paths $\lambda$ in $G$ based at the identity.  As for $P_eG^*$, $G$ acts on
  $\widehat{P_eG}$ by conjugation, and $\widehat{P_eG}$ is
  $G$-equivariantly contractible. Let $\widehat{\rho}: \widehat{P_eG}
  \to G$ denote the endpoint map; then we define the {\bf infinite
  approximating space $Y$} to be 
\begin{equation}\label{eq:infinite-approx-space}
Y := \left\{ (\lambda, m) \in \widehat{P_eG} \times M \, \vert \,
\Phi(m) = \widehat{\rho}(\gamma) \right\}. 
\end{equation}
As for \({\mathcal M}^*,\) projection to the first factor yields a
$G$-equivariant map \(\widehat{\mu}: Y \to \widehat{P_eG}.\)
Since on $\widehat{P_eG}$ the paths are piecewise smooth, we may also define on
$Y$ the following energy functional, familiar from classical Morse
theory, using the given bi-invariant metric on $G$: 
\begin{equation}\label{eq:def-hatf}
\widehat{f}(\lambda, m) := \int_0^1 \left\| \lambda^{-1} \frac{d\lambda}{dt}
\right\|^2 dt.
\end{equation}
Further, since
$\widehat{P_eG}$ is a subset of $P_eG^*$, there is a natural inclusion
\(Y \into {\mathcal M}^*,\) and again it is immediate from the definitions
that this inclusion takes $\widehat{f}^{-1}(0)$ exactly to
$(\mu^*)^{-1}({\mathbf e})$. 

To summarize, we have the following diagram of inclusions: 
\[
\xymatrix{
Y \ar @{^{(}->}[r] & {\mathcal M}^* & {\mathcal M} \ar @{_(->}[l] \\
\widehat{f}^{-1}(0) \ar @{^(->}[u] \ar[r]_{\cong} & (\mu^*)^{-1}(0)
\ar @{^(->}[u] &
\mu^{-1}(0) \ar @{^(->}[u] \ar[l]^{\cong}, \\
}
\]
where $Y, {\mathcal M}^*, {\mathcal M}$ are related by the diagram 
\[
\xymatrix{
Y \ar @{^{(}->}[r] \ar[d] & {\mathcal M}^* \ar[d] & {\mathcal M} \ar
@{_(->}[l] \ar[d] \\
M \ar[r]_{id} & M
 &
M  \ar[l]^{id}, \\
}
\]
which is obtained by pulling back the following diagram by the
quasi-Hamiltonian moment map \(\Phi: M \to G:\)
\[
\xymatrix{
\widehat{P_eG} \ar @{^{(}->}[r] \ar[d]_{\hat{\rho}} & {P_eG}^* \ar[d]_{\rho^*} & L_s\g^* \ar
@{_(->}[l] \ar[d]^{\Hol} \\
G \ar[r]_{id} & G
 &
G  \ar[l]^{id}. \\
}
\]
We may now state the main result of this section. 

\begin{proposition}\label{prop:surjection-approx-Y} 
Let $G$ be a compact Lie group equipped with a bi-invariant metric. Let $L_sG$ be the corresponding loop group for \(s > 3/2.\) Let \(({\mathcal M}, \omega, \mu)\) be  a Hamiltonian $L_sG$-space with proper moment map \(\mu: {\mathcal M} \to L_s\g^*.\) Let $Y$ denote the infinite
approximating space to ${\mathcal M}$ defined in~\eqref{eq:infinite-approx-space}, and
\(\hat{f}: Y \to \R\) the energy functional given in~\eqref{eq:def-hatf}. Then
${\mathcal M}$ and $Y$ are $G$-equivariantly homotopy-equivalent, and
\[
K^*_G({\mathcal M}) \to K^*_G(\mu^{-1}(0)) 
\]
is surjective if and only if 
\[
K^*_G(Y) \to K^*_G(\hat{f}^{-1}(0)) 
\]
is surjective. 
\end{proposition}

\begin{remark}
In \cite{BTW04}, Bott, Tolman, and Weitsman prove that the inclusions \(Y \into
{\mathcal M}^*\) and \({\mathcal M} \into {\mathcal M}^*\) induce
isomorphisms in $G$-equivariant cohomology \cite[Lemma
3.2]{BTW04}. Hence the content of our claim above is the stronger
statement that in fact these maps of spaces are
$G$-equivariant homotopy equivalences. 
\end{remark}

\begin{proof} 
%We now claim that \(\rho^*: P_eG^* \to G, \hat{\rho}: \widehat{P_eG}
%\to G,\) and \(\Hol: L_s\g^* \to G\) are all $G$-fibrations, where the
%action of $G$ on itself is by conjugation, and on each of the total
%spaces as described above. That they are fibrations in the
%non-equivariant sense, with fibers \(\Omega_eG^*,
%\widehat{\Omega_eG},\) and \(\Omega_sG\) respectively, is
%well-known. 

We prove the following sequence of geometric statements. 

\begin{enumerate} 
\item The maps \(\rho^*: P_eG^*
\to G,\) and \(\Hol: L_s\g^* \to G\) are both $G$-fibrations, where
the action of $G$ on itself is by conjugation, and the action on the
total spaces is by pointwise conjugation. Here, by a $G$-fibration we do not mean a principal $G$-bundle; rather, we mean a map of $G$-spaces satisfying a $G$-equivariant version of the standard homotopy lifting property (see for instance \cite{May96}). 

\item The total spaces of the pullback fibrations ${\mathcal M}^*$ and
${\mathcal M}$ are $G$-equivariantly
homotopy equivalent. 

\item The inclusion $\iota: \widehat{P_eG} \into P_eG^*$ is a
$G$-equivariant homotopy equivalence. Moreover, the $G$-equivariant homotopy inverse 
\(h: P_{e}G^{*} \to \widehat{P_{e}G}\) and the $G$-equivariant homotopies 
\(H^{*}: P_{e}G^{*} \times [0,1] \to P_{e}G^{*}\) and \(\widehat{H}: \widehat{P_{e}G} \times [0,1] 
\to \widehat{P_{e}G}\) with \(H^{*}(0) = \iota \circ h, H^{*}(1) = {\mathrm id}_{P_{e}G^{*}}, \widehat{H}(0) = h \circ \iota, \widehat{H}(1) = {\mathrm id}_{\widehat{P_{e}G}}\) may be chosen to preserve the projection maps to $G$, i.e.
\(\hat{\rho} \circ h = \rho^{*},\) and  \(\rho^{*} \circ H^{*}(t) = \rho^{*}, \widehat{\rho} \circ  \widehat{H}(t) = \widehat{\rho}\) for all \(t \in [0,1].\)

\item The total spaces of the pullback fibrations $Y$ and ${\mathcal
M}^*$ are $G$-equivariantly homotopy equivalent. 

\end{enumerate}

From these statements, it immediately follows that the spaces total
spaces $Y$ and ${\mathcal M}$ are $G$-equivariantly homotopy
equivalent. Moreover, since $G$-equivariantly homotopic spaces have
isomorphic equivariant $K$-theory, we have \(K^*_G({\mathcal M}) \cong
K^*_G({\mathcal M}^*) \cong K^*_G(Y).\) Furthermore, since we have
already seen that the $G$-equivariant inclusions \({\mathcal M} \into
{\mathcal M}^*\) and \(Y \into {\mathcal M}^*\) take \(\mu^{-1}(0)\)
exactly to \((\mu^*)^{-1}({\mathbf e})\) and \(\hat{f}^{-1}(0)\)
exactly to \((\mu^*)^{-1}({\mathbf e})\), respectively, we may also
immediately conclude that
\(K_G^*({\mathcal M}) \to K^*_G(\mu^{-1}(0))\) is surjective
if and only if \(K_G^*(Y) \to K^*_G(\hat{f}^{-1}(0))\) is
surjective. Hence in order to complete the proof, it suffices to prove
each of the geometric statements given above.

We begin with the first statement. That these are both fibrations in the
non-equivariant sense, with fibers \(\Omega_eG^*\) and \(\Omega_sG\) respectively, is
well-known. 
Here $\Omega_eG^*$ denotes the continuous loops based at the identity
and $\Omega_sG$ the maps \(S^1 \to G\) of Sobolev class $s$ and \(\gamma(0) = \gamma(1) = e \in G.\) 
That these are also
$G$-fibrations follows from the fact that here we take the $G$-action
to be trivial
on the domains $[0,1]$ and $S^1$ in the definition of the paths
\(\gamma: [0,1] \to G\) and loops \(\gamma: S^1 \to G.\) 
For the second statement, we first observe that both $L_s\g^*$ and
$P_eG^*$ are $G$-equivariantly homotopy equivalent to a point. In particular,
the inclusion \(\imath: L_s\g^* \into P_eG^*\) is a $G$-equivariant
homotopy equivalence, and by definition of the holonomy map, \(\Hol =
\rho^* \circ \imath.\) 
Now applying
Theorem~\ref{theorem:homotopic-fibrations-pullback} to the case where \(B=G, p = \Hol:
E = L_s\g^* \to G, p' = \rho^*: E' = P_eG^* \to G, \psi = \imath,
X = M\) the
quasi-Hamiltonian $G$-space associated to ${\mathcal M}$ and \(f =
\Phi\) the quasi-Hamiltonian moment map, we may immediately conclude
that ${\mathcal M}^*$ and ${\mathcal M}$ are $G$-equivariantly homotopy
equivalent. 

For the third statement, we construct explicitly the $G$-equivariant
homotopy inverse $h$ of the inclusion map \(\iota: \widehat{P_eG} \into
P_eG^*\) with the required properties. In the non-equivariant setting,
a
construction of an (ordinary) homotopy inverse to the 
inclusion \(\iota: \widehat{P_eG} \into P_eG^*\)
is given by Milnor in \cite[Section 17 and Appendix A]{Mil63}. Hence to prove
the statement, what remains to be seen is that the Milnor
argument can be made to be $G$-equivariant under the conjugation
action of $G$ on itself, in the sense stated above in (3). 
We will not reproduce
here the construction given by Milnor, as it is clearly written already in
\cite{Mil63}; we will instead only remark on those points which
require extra argument in our equivariant case. Following Milnor, we
call an open set $U$ in $G$ {\em geodesically convex} if any two
points $p,q$ in $U$ have a unique geodesic connecting them. We first
observe that any open covering ${\mathfrak U}$ of $G$ by geodesically
convex open sets can be assumed without loss of generality to be
$G$-invariant (under the conjugation action) since we may always add
to ${\mathfrak U}$ all $G$-translates of the open sets \(U \in
{\mathfrak U}.\) 
(Note that if $U$ is geodesically convex, then for any \(g \in G,\) so
is \(gUg^{-1},\) since the metric on $G$ is assumed bi-invariant.)
Following Milnor, we may now define $P_eG^*_k$ to be the subset
of $P_eG^*$ consisting of continuous paths $\gamma$ such that each subinterval
\([\frac{\ell-1}{2^k}, \frac{\ell}{2^k}], 1 \leq \ell \leq 2^k,\) is
mapped by $\gamma$ into one of the open sets in the cover ${\mathfrak
U}$. Since each  of the sets $U \in {\mathfrak U}$ 
is open in $G$, it is straightforward to see that under the 
uniform topology $P_eG^*_k$ is an open subset of $P_eG^*$. 
Moreover, since ${\mathfrak U}$ is $G$-invariant, it follows that
$P_eG^*_k$ is $G$-invariant, as is the
inverse image \(\widehat{P_eG}_k := \iota^{-1}(P_eG^*_k) \subseteq
\widehat{P_eG}.\) Again recalling that the metric in $G$ is
bi-invariant, it can now be verified that the homotopy inverse $h$ to
$\iota |_{\widehat{P_eG}_k}: \widehat{P_eG}_k \into P_eG^*_k$
constructed in
\cite[Section 17]{Mil63} is in fact also a $G$-equivariant homotopy
inverse. Moreover, it is evident from the definition of $h$ and the
argument in \cite[Section 16.2]{Mil63} that the homotopies to the
identity of \(h \circ \iota|_{\widehat{P_eG}_k}\) and
\(\iota|_{\widehat{P_eG}_k} \circ h\) are not only equivariant, but
also preserve the endpoint projection map for all \(t \in [0,1].\)

The next step is to check that, in our setting, the
argument in \cite[Example 1, Appendix A]{Mil63} that $\widehat{P_eG},
P_eG^*$ are $G$-homotopy direct limits of the $\widehat{P_eG}_k,
P_eG^*_k$ respectively, and that the proof of \cite[Theorem A, Appendix
A]{Mil63} can be made $G$-equivariant. Since $G$ is compact, we may
without loss of generality assume that the partition
of unity and $\R$-valued function $f$ in \cite[Example 1, Appendix
A]{Mil63} are $G$-invariant. Moreover, in the case that the original
space $X$ is equipped with a $G$-action and the sequence of subspaces
\(X_0 \subseteq X_1 \subseteq X_2 \cdots\) are all $G$-invariant (see \cite[p.169]{Mil63}) then
\(X \times \R\) can also be made a $G$-space with $G$ acting solely on
the first factor, and the space $X_{\Sigma}$ defined in
\cite[p.169]{Mil63} is $G$-invariant in \(X \times \R.\) The
projection \(p: X_{\Sigma} \to X\) then can be seen to be a 
$G$-equivariant homotopy equivalence, and we may conclude that 
$\widehat{P_eG},
P_eG^*$ are $G$-homotopy direct limits of the open $G$-invariant subsets $\widehat{P_eG}_k,
P_eG^*_k$ respectively, as desired. Similarly it is then
straightforward to check that a $G$-equivariant version of
\cite[Theorem A, Appendix A]{Mil63} is still valid, since with the
above assumptions the explicit constructions given in
\cite[p.150-153]{Mil63} are all $G$-equivariant. Moreover, it is 
evident from the construction that the $G$-equivariant homotopies preserve the
endpoint projection, as desired.

The fourth and last geometric statement follows from the following
more general statement. Suppose \(B, W, Z\) are $G$-spaces, and
\(\pi_W: W \to B\), \(\pi_Z: Z \to B,\) are $G$-equivariant maps. Suppose
also there exist $G$-equivariant maps \(\phi: W \to Z\), \(\psi: Z \to W\),
\(H_W: W\times [0,1] \to W\), \(H_Z: Z \times [0,1] \to Z\) (here the
$G$-action on the product manifolds are trivial on the second factor)
which satisfy \(\pi_Z \circ \phi = \pi_W\), \(\pi_Z \circ \psi = \pi_W,\)
and \(\pi_W \circ H_W(t) = \pi_W\), \(\pi_Z \circ H_Z(t) = \pi_Z\) for all
\(t \in [0,1].\) Now let \(\Phi^*W, \Phi^*Z\) denote the 
spaces obtained by pulling back via a map \(\Phi: M \to B\) the spaces
$W$ and $Z$, so 
\[
\Phi^*W := \{(m, w): \Phi(m) = \pi_W(w) \}, \quad 
\Phi^*Z := \{(m, z): \Phi(m) = \pi_Z(z) \}.
\]
Let the maps \(\Phi^*\phi:
\Phi^*W \to \Phi^*Z\), \(\Phi^*\psi: \Phi^*Z \to \Phi^*W\), \(\Phi^*H_W:
\Phi^*W \times [0,1] \to \Phi^*W\), \(\Phi^*H_Z: \Phi^*Z \times [0,1] \to
\Phi^*Z\) be defined respectively by 
\[
(m,w) \mapsto (m,\phi(w), \hs (m,z)
\mapsto (m, \psi(z)), \hs (m,w,t) \mapsto (m, H_W(w,t)), 
\hs (m,z,t) \mapsto
(m, H_Z(z,t)).
\]
Then it is straightforward to check that
these are well-defined, $G$-equivariant, and provide the necessary $G$-equivariant homotopy
equivalences and $G$-equivariant homotopies. 
Hence $\Phi^*W$ and $\Phi^*Z$ are also
$G$-equivariantly homotopy equivalent. Our case follows by taking \(W
= \widehat{P_eG}, Z = P_eG^*, B = G, \Phi:M \to G\) the
quasi-Hamiltonian moment map, and the $G$-equivariant homotopy
equivalences to be those constructed in the proof of the third
statement above.

\end{proof}

\section{Surjectivity via finite-dimensional
  approximations}\label{sec:finite}

In the previous section, we showed that surjectivity in $K$-theory for
the Hamiltonian $L_sG$-space ${\mathcal M}$ follows from an analogous
surjectivity result for the space $Y$. 
In this section, we argue that the surjectivity result for the space
$Y$ can in turn be reduced to that for an increasing sequence of
finite-dimensional approximating manifolds $Y_n$. Just as the space
$Y$ is constructed essentially by replacing the space of connections
$L_s\g^*$ with the space of piecewise smooth paths in $G$, we will
construct the finite-dimensional approximations $Y_n$ by replacing the
piecewise smooth paths in $G$ with piecewise geodesic paths in $G$ or
equivalently (if consecutive points are close enough), a sequence of points in $G$. 
We also define functions \(f_n: Y_n \to \R\) on each
finite-dimensional approximation $Y_n$ analogous to the energy
functional on paths. 

In order to then prove that the surjectivity
of~\eqref{eq:Kirwan-to-hatf-level} follows from the surjectivity of 
\begin{equation}\label{eq:Yn-to-fn-level} 
K^*_G(Y_n) \to K^*_G(f_n^{-1}(0))
\end{equation}
for each $n$, where the maps are induced by the inclusions
\(f_n^{-1}(0) \into Y_n,\) we will show that in fact each $Y_n$ is
$G$-equivariantly homotopy-equivalent to the preimage
$\hat{f}^{-1}(-\infty, cn) \subseteq Y$ for a positive constant $c$. We also show
that the appropriate 
Milnor $\lim^1$ term vanishes when comparing the limit of the
equivariant $K$-theories of the $Y_n$ to that of $Y$. This vanishing
turns out to be a straightforward consequence of the geometric
arguments given in \cite{BTW04} and the $K$-theoretic Atiyah-Bott
lemma \cite[Lemma 2.3]{HarLan07}; moreover, it allows us to quickly obtain the
necessary final step in the argument, namely,
that~\eqref{eq:Yn-to-fn-level} is indeed surjective for each $Y_n$. 

As was the case for $Y$, our construction of the $Y_n$ is the same as
that given in \cite[Section 4]{BTW04} so we only briefly recall it here and refer
the reader to the original work for details; however, 
the argument that the relevant ring surjections are related does 
require extra argument and care in the $K$-theory case. 
For \(n
\in \N\) we define 
\[
X_n := \left\{ (g_1, g_2, \ldots, g_n, m) \in G^n \times M \,
\big\vert \, \Phi(m) = g_n \right\} \cong G^{n-1} \times M, 
\]
where $M$ is our quasi-Hamiltonian $G$-space with moment map
$\Phi$. The space $X_n$ is naturally a $G$-space, given by the
diagonal action; here $G$ acts on itself by conjugation. We define the
analogue of the energy function \(f_n: X_n \to \R\) by 
\[
f_n(g_1, g_2, \ldots, g_n, m) := n\rho(e,g_1)^2 + n\rho(g_1, g_2)^2 +
\cdots + n \rho(g_{n-1}, g_n)^2,
\]
where $\rho$ denotes the distance between two points on $G$ with
respect to the given metric. Let \(\bar{\rho} > 0\) be a positive real
number such that if two points $p, q \in G$ are of distance less than
$\bar{\rho}$ from each other, i.e. \(\rho(p,q) < \bar{\rho},\) then 
there exists a unique shortest geodesic in $G$ 
connecting them. We then define the {\bf finite-dimensional
  approximating manifolds} to $Y$ by 
\begin{equation}\label{eq:finite-approximating-spaces}
Y_n := f_n^{-1}\left(-\infty, \frac{1}{2}n\bar{\rho}^2\right).
\end{equation}
Since $f_n$ is $G$-invariant, $Y_n$ is also a $G$-space.

Our goal in this section is to analyze the infinite-dimensional space
$Y$ using the sequence of finite-dimensional approximations
$Y_n$. However, in order to effectively accomplish this, we must
understand the relationship between the successive approximations. As
a first step, 
we briefly recall the
relationship, explained in detail in \cite[Section 4]{BTW04}, between
$Y_n$ and \(\hat{f}^{-1}(-\infty, \frac{1}{2}n\bar{\rho}^2).\) The
essential observation is that if 
\[
f_n(g_1, g_2, \ldots, g_n, m) = n \left( \sum_{i=0}^{n-1} \rho(g_i,
g_{i+1})^2 \right) < \frac{1}{2} n\bar{\rho}^2,
\]
where we set \(g_0 = e \in G,\) then \(\rho(g_i, g_{i+1}) <
\bar{\rho}\) for all \(0 \leq i \leq n-1.\) In other words, in this situation the
consecutive group elements have a unique shortest
geodesic connecting them. By joining these geodesics together (in $n$
equal subintervals) to form a
single path from $e$ to $g_n$ we may define a map \(\beta_n:
Y_n \to Y.\) Since each geodesic from time $[\frac{i}{n},
\frac{i+1}{n}]$ has constant speed \(\frac{\rho(g_k, g_{k+1})}{1/n} =
n\rho(g_k, g_{k+1}),\) it is straightforward to see that the path from
$[0,1]$ has total energy \(n\left(\sum_{i=0}^{n-1}
\rho(g_i,g_{i+1})^2\right)\) and hence 
the map $\beta_n$ satisfies \(\hat{f} \circ \beta_n = f_n.\) So in fact
$\beta_n$ maps $Y_n$ into $\hat{f}^{-1}(-\infty,
\frac{1}{2}n\bar{\rho}^2) \subset Y$. A reverse map \(\alpha_n:
\hat{f}^{-1}(-\infty, \frac{1}{2}n\bar{\rho}^2) \to Y_n\) can be
defined by taking an element \((\lambda, m) \in Y := \left\{ (\lambda,
m) \in 
\widehat{P_eG} \times M \, \big\vert \, \lambda(1) = \Phi(m) \in G
\right\}\) to the element \((\lambda(\frac{1}{n}),
\lambda(\frac{2}{n}), \ldots, \lambda(\frac{n-1}{n}), \lambda(1), m)
\in Y_n.\) It is shown in \cite[Lemma 4.2]{BTW04} that $\alpha_n$ is
well-defined, and that \(\beta_n: Y_n \to \hat{f}^{-1}(-\infty,
\frac{1}{2}n\bar{\rho}^2)\) is a $G$-homotopy equivalence with
homotopy inverse $\alpha_n$. In particular, we may conclude
that
\begin{equation}\label{eq:Y_n-preimage}
K^*_G\left(\hat{f}^{-1}\left(-\infty, \frac{1}{2}n\bar{\rho}^2\right)\right) \cong
K^*_G(Y_n). 
\end{equation}
Furthermore, using the natural inclusion \(\hat{f}^{-1}\left(-\infty,
\frac{1}{2}n\bar{\rho}^2\right) \into \hat{f}^{-1}\left(-\infty,
\frac{1}{2}(n+1)\bar{\rho}^2\right),\) we may now define a map
\(\phi_n: Y_n \to Y_{n+1}\) between the successive approximations as
the composition 
\begin{equation}\label{eq:def-phi_n} 
\xymatrix{
Y_n \ar[r]^(0.25){\beta_n} & \hat{f}^{-1}\left(-\infty, \frac{1}{2}n\bar{\rho}^2\right)
\ar @{^{(}->}[r] & \hat{f}^{-1}(-\infty, \frac{1}{2}(n+1)\bar{\rho}^2)
\ar[r]^(0.70){\alpha_{n+1}} & Y_{n+1}.
}
\end{equation}

We may now state our first result of the section, which relates the equivariant
$K$-theories of the successive approximations. 

\begin{proposition}\label{prop:surjection-for-each-n} 
For each $n \in \N$, the map \(\phi_n: Y_n \to Y_{n+1}\)
in~\eqref{eq:def-phi_n} induces a surjection in
equivariant $K$-theory, i.e. the induced ring homomorphism 
\[
\phi_n^*: K^*_G(Y_{n+1}) \to K^*_G(Y_n)
\]
is surjective. 
\end{proposition}

The proof, as we see below, is by using the Morse-Kirwan theory of the
function $f_{n+1}$ on $Y_{n+1}$. 

\begin{proof} 
We will first prove the following
geometric statements: 
\begin{enumerate} 
\item \(\phi_n: Y_n \into Y_{n+1}\) is an inclusion; 
\item \(\phi_n(Y_n) \subseteq f_{n+1}^{-1}(-\infty,
  \frac{1}{2}n\bar{\rho}^2)\); and 
\item \(\phi_n: Y_n
  \to f^{-1}_{n+1}\left(-\infty, \frac{1}{2}n\bar{\rho}^2\right)\) 
is a $G$-equivariant homotopy equivalence. 
\end{enumerate}

We begin with the first statement.  Let \((g_1, g_2, \ldots,
g_n, m)\) and \((g'_1, g'_2, \ldots, g'_n, m')\) be two distinct
elements in $Y_n$. It is evident from the definitions that $\phi_n$ does not change the last
coordinate.  Hence if
\(m \neq m'\) then it is immediate that \(\phi(g_1, g_2, \ldots, g_n,
m) \neq \phi(g'_1, g'_2, \ldots, g'_n, m'),\) so we may suppose that
there exists an index $k$, with \(1 \leq k \leq n,\) such that \(g_k \neq g'_k.\) In fact let $k$
denote the least such index. By definition of $\phi_n$, the $k$-th
component of $\phi_n(g_1, \ldots, g_n, m)$ is given by
\(\lambda(\frac{k}{n+1}),\) where \(\lambda: [0,1] \to G\) is the
based path starting at \(e \in G, \lambda(\frac{i}{n}) = g_i\) for \(1
\leq i \leq n,\) and \(\lambda: \left[ \frac{i-1}{n}, \frac{i}{n}
  \right] \to G\) is the unique geodesic connecting $g_{i-1}$ and
$g_i$ (we set \(g_0 = e \in G\)). Let $\lambda'$ denote the image
\(\phi_n(g'_1, g'_2, \ldots, g'_n, m').\) By assumption on the index
$k$, \(g_{i} = g'_{i}\) for \(i \leq k-1.\) Since
\(\rho(g_{k-1}, g_k) < \bar{\rho}\) and \(\rho(g'_{k-1}, g'_k) <
\bar{\rho},\) there exists a unique shortest geodesic \(g_{k-1}\cdot
\exp(tX)\) connecting $g_{k-1}$ and $g_k$, and similarly
\(g_{k-1}\exp(tX') = g'_{k-1}\exp(tX')\) connecting \(g_{k-1} =
g'_{k-1}\) and $g'_k$, both obtained by a translate (by multiplication
on $G$) of a $1$-parameter subgroup in $G$. Since \(g_k \neq g'_k,\)
we must have \(X \neq X',\) and the geodesics cannot intersect. Hence
the image of \(\lambda\left(\left[\frac{k-1}{n},
  \frac{k}{n}\right]\right)\) is disjoint from
\(\lambda'\left(\left[\frac{k-1}{n}, \frac{k}{n}\right]\right),\) and
in particular \(\lambda\left(\frac{k}{n+1}\right) \neq
\lambda'\left(\frac{k}{n+1}\right).\) Hence $\phi_n$ is injective. 

For the second statement, to see that \(\phi_n(Y_n)\) is a subset of \(f_{n+1}^{-1}(-\infty,
\frac{1}{2}n\bar{\rho}^2),\) it suffices to recall that the inclusion
\(\beta_n: Y_n \into \hat{f}^{-1}\left(-\infty, \frac{1}{2}n\bar{\rho}^2\right)\)
satisfies \(\hat{f} \circ \beta_n = f_n,\) and that the retraction
\(\alpha_{n+1}: \hat{f}^{-1}(-\infty, \frac{1}{2}(n+1)\bar{\rho}^2)
\to Y_{n+1}\) does not increase energy, since a geodesic between
\(\lambda\left(\frac{i-1}{n+1}\right)\) and
\(\lambda\left(\frac{i}{n+1}\right)\) always has energy less than or
equal to that of any path \(\lambda: \left[ \frac{i-1}{n+1},
  \frac{i}{n+1}\right] \to G.\) Hence if \((g_1, g_2, \ldots, g_n, m)
\in Y_n,\) i.e. \(f_n(g_1, \ldots, g_n, m) \leq
\frac{1}{2}n\bar{\rho}^2,\) then \(f_{n+1} \circ \phi_n(g_1, \ldots,
g_n, m)\) is also \(\leq \frac{1}{2}n\bar{\rho}^2,\) as desired. 

Finally, we prove that 
$\phi_n$ is a $G$-equivariantly homotopy equivalence between $Y_n$ and the subset 
\(f_{n+1}^{-1}\left(-\infty, \frac{1}{2}n\bar{\rho}^2\right) \subseteq
Y_{n+1}.\)
Recall that 
the inclusion \(\beta_n: Y_n \to \hat{f}^{-1}(-\infty,
\frac{1}{2}n\bar{\rho}^2)\) is already known to be a $G$-homotopy
equivalence, as described above. 
Furthermore, 
since the
$G$-equivariant retraction $\alpha_{n+1}$ does not increase energy, as noted
above, we also have a well-defined restriction 
\begin{equation}\label{eq:res-alpha}
\alpha_{n+1}\vert_{\hat{f}^{-1}\left(-\infty, \frac{1}{2}n\bar{\rho}^2\right)}: \hat{f}^{-1}\left(-\infty, \frac{1}{2}n\bar{\rho}^2\right) \to
f_{n+1}^{-1}\left(-\infty, \frac{1}{2}n\bar{\rho}^2\right).
\end{equation}
Putting these together, we see that it suffices to show that this
map~\eqref{eq:res-alpha} is a $G$-equivariant homotopy equivalence.

By construction, we have that \(\alpha_{n+1} \circ \beta_{n+1} =
\mathrm{id},\) so it suffices to show that the restriction
\(\beta_{n+1} \circ \alpha_{n+1}\) to $\hat{f}^{-1}(-\infty,
\frac{1}{2}n\bar{\rho}^2)$ is $G$-equivariantly homotopic to the
identity. In \cite[Lemma 4.2]{BTW04}, a $G$-equivariant homotopy $D_t$
between \(\beta_{n+1} \circ \alpha_{n+1}\) and the identity map on
\(\hat{f}^{-1}(-\infty, \frac{1}{2}(n+1)\bar{\rho}^2)\) is explicitly
constructed as follows. 
Let \(t \in [0,1]\) and \(\lambda \in \hat{f}^{-1}\left(-\infty,
\frac{1}{2}(n+1)\bar{\rho}^2\right).\) Then $D_t(\lambda): [0,1] \to
G$ is defined as follows: for every index $i$, \(0 \leq i \leq n,\) on
the interval \([\frac{i}{n+1}, \frac{i+t}{n+1}],\) $D_t(\lambda)$ is
the unique geodesic joining \(\lambda(\frac{i}{n+1})\) and
\(\lambda(\frac{i+t}{n+1})\). On the interval \([\frac{i+t}{n+1},
\frac{i+1}{n+1}],\) $D_t(\lambda)$ is defined to be the same as the
original path \(\lambda: [\frac{i+t}{n+1}, \frac{i+1}{n+1}] \to G.\)
It is evident that at \(t=0,\) \(D_0 \equiv id,\) and \(D_1 =
\beta_{n+1} \circ \alpha_{n+1}.\) Moreover, $D_t$ cannot increase the
energy of a path for any \(t \in [0,1],\) 
so for each $t$ the map $D_t$ preserves the subset
$\hat{f}^{-1}(-\infty,
\frac{1}{2}n\bar{\rho}^2)$.		In particular, it provides a 
$G$-equivariant homotopy between the restriction of $\beta_{n+1} \circ
\alpha_{n+1}$ to $\hat{f}^{-1}\left(-\infty, \frac{1}{2}n\bar{\rho}^2\right)$ and
the identity map, as desired. 

To conclude the proof, we must now show that the inclusion \(\phi_n:
Y_n \to Y_{n+1}\) induces a surjection in equivariant $K$-theory. From
the arguments above, it in fact suffices to show that the inclusion
\(f_{n+1}^{-1}(-\infty, \frac{1}{2}n\bar{\rho}^2) \into Y_{n+1}\)
induces a surjection in equivariant $K$-theory. We
first recall that 
$f_{n+1}$ is a Morse-Kirwan function, as is shown in \cite[Sections
  5-8]{BTW04}. This fact is highly non-trivial, but as it is explained
in detail in the above-mentioned reference, we will not discuss it
further here. We only note in particular that \cite[Proposition
  8.2]{BTW04} additionally proves the following geometric fact. Let
$C$ be a connected component of the critical set of $f_{n+1}$, and
$E_C^{-}$ be its negative normal bundle with respect to
$f_{n+1}$. Then they show that there exists a subtorus $T$ of $G$ and
a $Z(T)$-invariant subset $B$ of $C^T$ such that the natural map \(G
\times_{Z(T)} B\) is a $G$-equivariant homeomorphism; moreover,
$(E_C^{-})^T$ is a subset of the zero section of $E_C^-$. In this
geometric situation, we may immediately apply the $K$-theoretic
Atiyah-Bott lemma \cite[Lemma 2.3]{HarLan07} and conclude that the
$G$-equivariant $K$-theoretic Euler class $e_G(E_C^-) \in K^*_G(C)$ is
not a zero divisor for all components $C$. 
We briefly recall the geometric implications of this last statement about the Euler class, referring the reader to \cite{BTW04, HarLan07} for details. For sufficiently small \(\varepsilon > 0\) let $Y_{n+1}^{\pm}$ denote $f_{n+1}^{-1}(-\infty, f_{n+1}(C) \pm \varepsilon)$. There is then a commutative diagram involving the long exact sequence in equivariant $K$-theory of the pair $(Y_{n+1}^+, Y_{n+1}^-)$ as follows: 
\[
\xymatrix{ 
\cdots \ar[r] & K^*_G(Y_{n+1}^+, Y_{n+1}^-) \ar[r] \ar[d]_{\cong} & K^*_G(Y_{n+1}^+) \ar[r] \ar[d] & K^*_G(Y_{n+1}^-) \ar[r] & \cdots \\
 & K^{*-\lambda_C}_G(C) \ar[r]_{\cup e_G(E_C^-)} & K^*_G(C) & & \\
}
\]
where $\lambda_C$ denotes the Morse index of $C$. Since $e_G(E_C^-)$ is not a zero divisor, the bottom horizontal arrow is an injection, which in turn implies that the top long exact sequence splits. In particular, the restriction map \(K^*_G(Y^+_{n+1}) \to K^*_G(Y_{n+1}^-)\) is a surjection. Since this is true for all components $C$ of the critical set $\Crit(f_{n+1})$, we conclude that the restriction map \(K^*_G(Y_{n+1}) \to K^*_G(f_{n+1}^{-1}(-\infty, a))\) for any value \(a \in \R\) is a surjection. Taking \(a = \frac{1}{2}n\bar{\rho}^2,\) we then see that the composite ring homomorphism 
\[
\xymatrix{
K^*_G(Y_{n+1}) \ar[r] &  K^*_G(f_{n+1}^{-1}(-\infty,
\frac{1}{2}n\bar{\rho}^2)) \ar[r]^(0.65){\cong}_(0.65){\phi^*} &  K^*_G(Y_n)
}
\]
is a surjection, as desired. 
\end{proof}

We may now compare the limit of the $Y_n$ with all of $Y$, and also
conclude that the ring map~\eqref{eq:Kirwan-to-hatf-level} is a surjection if the
analogous map is a surjection for all $n$.  

\begin{proposition}\label{prop:finite-surjectivity} 
Suppose that for each \(n \in \N,\)
the inclusion $f_n^{-1}(0) \into Y_n$ induces a ring surjection 
\[
K^*_G(Y_n) \onto K^*_G(f_n^{-1}(0)).
\]
Then the natural ring homomorphism induced by the inclusion
\(\hat{f}^{-1}(0) \into Y,\)
\[
K^*_G(Y) \to K^*_G(\hat{f}^{-1}(0)),
\]
is also a surjection. 
\end{proposition}

\begin{proof} 
The union of the infinite increasing sequence of
$G$-invariant subspaces
\[
\cdots \subset \hat{f}^{-1}(-\infty, \frac{1}{2}(n-1)\bar{\rho}^2)
\subset \hat{f}^{-1}(-\infty, \frac{1}{2}n\bar{\rho}^2) \subset
\hat{f}^{-1}(-\infty, \frac{1}{2}(n+1)\bar{\rho}^2) \subset \cdots
\]
is all of $Y$. The $G$-equivariant inclusions then induce a sequence of
ring homomorphisms 
\[
\cdots \to K^*_G(\hat{f}^{-1}(-\infty, \frac{1}{2}(n+1)\bar{\rho}^2))
\to K^*_G(\hat{f}^{-1}(-\infty, \frac{1}{2}n\bar{\rho}^2)) \to
K^*_G(\hat{f}^{-1}(-\infty, \frac{1}{2}(n-1)\bar{\rho}^2)) \to \cdots 
\]
with inverse limit \(\underleftarrow{\lim} \hsm K^*_G(\hat{f}^{-1}(-\infty,
\frac{1}{2}n\bar{\rho}^2)).\) 
We now wish 
to compare this inverse limit with \(K^*_G(Y)\) 
using the Milnor sequence (see e.g. \cite[Theorem 13.1.3]{Sel97}). Note that $K^*_G$ satisfies the wedge axiom since by definition 
it is a representable theory; moreover, each of the subsets $\hat{f}^{-1}(-\infty, \frac{1}{2}n\bar{\rho}^2)$ 
are open subsets so we may use the Mayer-Vietoris sequence as in the proof of \cite[Theorem 13.1.3]{Sel97}. Hence we may conclude that there exists a sequence (``the Milnor sequence'') 
\begin{equation}\label{eq:Milnor} 
0 \to \underleftarrow{\lim}^1 K^*_G ( \hat{f}^{-1}(-\infty, \frac{1}{2}n\bar{\rho}^2)) 
\to K^*_G(Y) \to \underleftarrow{\lim} K^*_G(\hat{f}^{-1}(-\infty,
\frac{1}{2}n\bar{\rho}^2)) \to 0.
\end{equation}
We have just seen in Proposition~\ref{prop:finite-surjectivity} that for each
$n$, the map 
\[
K^*_G(Y_{n+1}) \onto K^*_G(Y_n)
\]
is surjective.  Since we have also just seen in the proof of
Proposition~\ref{prop:finite-surjectivity} that each $Y_n$ is
$G$-equivariantly homotopic to $\hat{f}^{-1}(-\infty,
\frac{1}{2}n\bar{\rho}^2)$, we conclude that 
\[
K^*_G( \hat{f}^{-1}(-\infty, \frac{1}{2}(n+1)\bar{\rho}^2)) \onto
K^*_G( \hat{f}^{-1}(-\infty, \frac{1}{2}n\bar{\rho}^2))
\]
is also surjective for each $n$. Hence the Mittag-Leffler condition is satisfied and the Milnor $\lim^1$ term, i.e. the first term in the short exact sequence~\eqref{eq:Milnor}, vanishes and 
we conclude
\begin{equation}\label{eq:inverse-lim}
K^*_G(Y) \cong 
\underleftarrow{\lim} K^*_G(\hat{f}^{-1}(-\infty,
\frac{1}{2}n\bar{\rho}^2)).
\end{equation}
Finally, we observe that each $f_n^{-1}(0)$ for all \(n \in \N\) is
identified via $\beta_n$ with $\hat{f}^{-1}(0)$, and that the maps
$\phi_n$ also identify each $f_n^{-1}(0) \subseteq Y_n$ with
\(f_{n+1}^{-1}(0) \subseteq Y_{n+1}.\) The result now follows from the
definition of the inverse limit. 
\end{proof}

We may now prove Theorem~\ref{theorem:main} as a straightforward corollary.

\begin{proof} ({\bf Proof of Theorem~\ref{theorem:main}}) 
By Propositions~\ref{prop:surjection-approx-Y},
\ref{prop:surjection-for-each-n}, and \ref{prop:finite-surjectivity},
it suffices to prove that for each \(n \in \N,\) the inclusion
$f_n^{-1}(0) \into Y_n$ induces a surjection in equivariant
$K$-theory.  The proof of Proposition~\ref{prop:finite-surjectivity} shows that $f_n$ is
a self-perfecting Morse-Kirwan function in equivariant
$K$-theory. Hence, in particular, since $f_n^{-1}(0)$ is the minimum
level of $f_n$, the inclusion \(f_n^{-1}(0) \into Y_n\) indeed induces
a surjection
\[
K^*_G(Y_n) \onto K^*_G(f_n^{-1}(0))
\]
in equivariant $K$-theory, as desired. 
\end{proof}

\appendix
\section{Some $G$-equivariant homotopy theory} 

In this Appendix we present, in as streamlined a fashion as possible, a small amount of 
$G$-equivariant homotopy theory which we require in order to make
our arguments in Section~\ref{sec:infinite}. We suspect that
such arguments as we record below are standard, well-known, and even
trivial to the experts,
but we were unable to find complete proofs. As such, with an apology
to those experts, we include them here. 

Our main references for the non-equivariant theory are the monographs by
Selick \cite{Sel97} and May \cite{May75}. All of the results we state below are obtained by
carefully inserting the $G$-action into the relevant statements in
\cite[Chapter 7]{Sel97}. We begin with some general facts about
$G$-fibrations and $G$-cofibrations, such as the $G$-equivariant version of a standard but crucial
fact, namely, that any $G$-map can be factored as a composition of a
$G$-fibration and a $G$-homotopy
equivalence. The main result for our purposes is
Theorem~\ref{theorem:homotopic-fibrations-pullback}, which states that if two $G$-fibrations are
$G$-homotopy equivalent, then their pullbacks along any $G$-map are
also $G$-homotopy equivalent.

\begin{theorem}\label{theorem:factor-fibration}
Let \(f: X \to Y\) be a based $G$-map which induces a surjective map
on path components. Then there exists a factorization \(f = p\phi\)
where \(\phi: X \to P_f\) is a $G$-homotopy equivalence and \(p: P_f
\to Y\) is a $G$-fibration.
\end{theorem}

The proof of this theorem follows the ``standard'' proof, which
explicitly constructs the intermediate space $P_f$ and maps \(\phi: X
\to P_f, p: P_f \to Y\) that satisfies the required properties. The
only extra work is to take care to check the $G$-equivariance
properties of this standard construction; this is what we do below.

\begin{proof} 
The inclusion \(\{0\} \cup \{1\} \into I = [0,1]\) is a cofibration in
the non-equivariant sense.
By equipping both the domain and the target with the
trivial $G$-action, the inclusion becomes a $G$-cofibration. It follows from the $G$-equivariant exponential law \cite[p.11]{May96}
that the induced map \(\ev_0 \times \ev_1: \Map(I,Y) \to \Map(\{0\}
\cup \{1\}, Y) = Y \times Y\) given by evaluation at $0$ and $1$,
respectively, is a $G$-fibration. Hence it follows that \(1_X
\times \ev_0 \times \ev_1: X \times \Map(I,Y) \to X \times Y \times
Y\) is also a $G$-fibration. Now define the map \(q: X \times Y \to X
\times Y \times Y\) by \(q(x,y) := (x, f(x), y).\) Since $f$ is by
definition $G$-equivariant, this is a $G$-equivariant map, and hence
we obtain a $G$-fibration $P_f$ over \(X \times Y\) by pulling back
\(1_X \times \ev_0 \times \ev_1\) by $q$.  More specifically, we have
\begin{align*}
P_f & := \left\{ (x,y, x', \gamma) \in X \times Y \times X \times \Map(I,Y) \, \bigg\vert \, x' = x, \gamma(0) = f(x), \gamma(1) = y \right\} \\
 & \cong \left\{ (x,\gamma) \, \bigg\vert \, \gamma(0) = f(x) \right\},\\
\end{align*}
with projection map \(\bar{p}: P_f \to X \times Y\) given by \(\bar{p}(x,\gamma) = \gamma(1).\) Since the map \(X \times Y \to Y\) which projects to the second factor is also a $G$-fibration, we obtain the map \(P_f \to Y\) by composing $\bar{p}$ with this projection, i.e. \(p: P_f \to Y\) is defined by \(p(x,\gamma) = \gamma(1)\) and is a $G$-fibration, since compositions of $G$-fibrations is also a $G$-fibration. 

It now suffices to prove that there exists a map \(\phi: X \to P_f\)
which is a $G$-homotopy equivalence. We first define \(\phi: X \to
P_f\) by \(\phi(x) := (x, c_{f(x)}),\) where $c_{f(x)}$ denotes the
constant path at \(f(x) \in Y.\) Since we have equipped the interval
$I$ with the trivial $G$-action, $\phi$ is also $G$-equivariant. Now
define \(\psi: P_f \to X\) by \(\psi(x,\gamma) = x.\) This is clearly
$G$-equivariant. We have by definition that \(\psi \circ \phi = 1_X,\) so to obtain
the $G$-homotopy-equivalence, it suffices to show that \(\phi \circ
\psi\) is $G$-homotopic to the identity $1_{P_f}$. We may explicitly
construct such a $G$-homotopy \(H: P_f \times I \to P_f\) by defining
\(H(x,\gamma,s) := (x, \gamma_s), \) where for any path
$\gamma$ the adjusted path $\gamma_s$ is given by \(\gamma_s(t) :=
\gamma(st).\) It is now straightforward to check that $H$ is a
$G$-equivariant homotopy between \(\phi \circ \psi\) and $1_{P_f}$, as desired.

\end{proof}

\begin{theorem}\label{theorem:homotopy-pullback} 
Let \(p: E \to B\) be a $G$-fibration, and let \(f: X \to B, g: X \to
B\) be $G$-maps which are $G$-homotopy equivalent. Then the pullback
fibrations \(f^*E \to X, g^*E \to X\) are $G$-fibre-homotopy equivalent. 
\end{theorem}

\begin{proof} 
We will explicitly construct the $G$-fibre-homotopy
equivalence. Throughout, we will use $\mathrm{pr}_i$ to denote the
projection of a direct product to its $i$-th factor. First, by writing
$f$ and $g$ as composites with the $G$-homotopy, we may without loss
of generality consider the special case in which the base $B$ is of
the form \(X \times I\) and \(f,g\) are the inclusions at the two ends
\(\iota_0, \iota_1.\) Here $X$ is a $G$-space, $I$ is equipped with
the trivial $G$-action, and \(X \times I\) the diagonal
$G$-action. Let \(E_s := p^{-1}(X \times \{s\})\) denote the pullback
of $p$ by \(\iota_s: X \to X \times \{s\}.\) Let \(h: X \times I
\times I \times I \to X \times I\) be defined by \(h(x,r,s,t) = (x,
(1-t)r+st).\) All factors of $I$ are equipped with the trivial
$G$-action so this is clearly $G$-equivariant. 

Since \(p: E \to X \times I = B\) is a $G$-fibration by assumption, we
may apply the homotopy lifting property to the following diagram: 
\[
\xymatrix{
E \times I \ar[rr]^{\mathrm{pr}_1}  \ar[d]^{\iota_{E,0}} && E \ar[d]^p
\\
E \times I \times I \ar[rr]_{h \circ (p \times 1 \times
1)} && X \times I 
}
\]
where \(\iota_{E,0}: E \times I \to E \times I \times I\) is the
inclusion \((e,s) \mapsto (e,s,0)\) and \(p \times 1 \times
1(e,s,t) = (p(e), s, t) \in (X \times I) \times I \times
I.\) As a result, we obtain a diagonal map \(F: E \times I \times I
\to E\) making both resulting triangles commute. We set \(K :=
F(-,-,1): E \times I \to E\) given by \(K(e,s) = F(e,s,1).\) Observe
that the map \(H: E \times I \to E\) given by \(H(e,t) = F(e,
\mathrm{pr}_2 \circ p(e), t)\) gives a $G$-homotopy between
\(1_E\) (at \(t = 0\)) and the $G$-map \(k(e) = K(e,
\mathrm{pr}_2 \circ p(e)).\) Moreover, by construction it satisfies
\(pH(e,t) = p(e) \forall t.\) Hence, in order to show that $E_0$ and
$E_1$ are $G$-fibre-homotopy-equivalent, it suffices to construct
\(\alpha: E_0 \to E_1, \beta: E_1 \to E_0\) such that \(\beta \circ
\alpha\) and \(\alpha \circ \beta\) are $G$-homotopic to \(k \circ k
|_{E_0}\) and \(k \circ k|_{E_1}\) respectively. 

We define \(\alpha: E_0 \to E_1\) by \(\alpha(e) = K(e,1)\) and
\(\beta: E_1 \to E_0\) by \(\beta(e) = K(e,0);\) these are
well-defined since \(pK(e,s) = (\mathrm{pr}_1 \circ p(e), s)\) so in
particular \(K(e,s) \in E_s.\) Then it is straightforward to check
that \((e,s) \mapsto K(K(e,1-s),0)\) and \((e,s) \mapsto K(K(e,s),1)\)
provide $G$-homotopies from \(\beta \circ \alpha\) to \(k \circ
k|_{E_0}\) and \(\alpha \circ \beta\) to \(k \circ k|_{E_1}\)
respectively, both covering the constant homotopy \(1_X
\sim 1_X.\) The claim follows. 
\end{proof}

Given this theorem, we can easily prove the following. It is obtained
by transforming an arbitrary $G$-map into a $G$-fibration by
Theorem~\ref{theorem:factor-fibration} and then applying
Theorem~\ref{theorem:homotopy-pullback} three times.

\begin{theorem}\label{theorem:homotopic-fibrations-pullback}
Let \(p: E \to B\) and \(p': E' \to B\) be $G$-fibrations. Suppose
there exists a $G$-homotopy equivalence \(\psi: E \to E'\) satisfying
\(p' \circ \psi = p.\) Then for any $G$-equivariant map \(f: X \to
B,\) the pullback $G$-fibrations \(f^*E \to X, f^*E' \to X\) are also
$G$-homotopy equivalent. 
\end{theorem}

\begin{proof} 
By Theorem~\ref{theorem:factor-fibration}, we may decompose $f$ into a
composition \(f = f' \circ \phi,\) 
\[
\xymatrix{
X \ar[r]^{\phi} & P_f \ar[r]^{f'} & B \\
}
\]
where $f'$ is a $G$-fibration and $\phi$ is a $G$-homotopy
equivalence. Then denote by $Q_f$ and $Q'_f$ the pullbacks of $E$ and
$E'$ via the map $f'$, i.e. we have the fiber squares
\[
\xymatrix{
Q_f \ar[r] \ar[d] & E \ar[d]^p \\
P_f \ar[r]_{f'} & B 
}
\quad \mbox{and} \quad 
\xymatrix{
Q'_f \ar[r] \ar[d] & E' \ar[d]^{p'} \\
P_f \ar[r]_{f'} & B. 
}
\]
Since \(f': P_f \to B\) is a $G$-fibration, and \(p: E \to B\) and
\(p': E' \to B\) are $G$-homotopy equivalent, we may now apply
Theorem~\ref{theorem:homotopy-pullback} to conclude that $Q_f$ and
$Q'_f$ are also $G$-homotopy equivalent. 
Now consider the pullbacks $Q$ and $Q'$ of $E$ and $E'$, respectively,
by $f$, i.e. we have the fiber squares
\[
\xymatrix{
Q \ar[r] \ar[d] & E \ar[d]^p \\
X \ar[r]_{f} & B 
}
\quad \mbox{and} \quad 
\xymatrix{
Q' \ar[r] \ar[d] & E' \ar[d]^{p'} \\
X \ar[r]_{f} & B. 
}
\]
Since \(f = f' \circ \phi,\) these fiber squares in fact also fit into
the larger commutative diagrams of pullbacks 
\begin{equation}\label{eq:full-pullback} 
\xymatrix{
Q \ar[r] \ar[d] & Q_f \ar[r] \ar[d] & E \ar[d]^p \\
X \ar[r]_{\phi} & P_f \ar[r]_{f'} & B 
}
\quad \mbox{and} \quad 
\xymatrix{
Q' \ar[r] \ar[d] & Q'_f \ar[r] \ar[d] & E' \ar[d]^{p'} \\
X \ar[r]_{\phi} & P_f \ar[r]_{f'} & B. 
}
\end{equation}
Looking only at the left square of each of these diagrams, we have
that the base $\phi$ is a $G$-homotopy equivalence by assumption. A
simple corollary of Theorem~\ref{theorem:homotopy-pullback} is that if
a map \(h: A \to B\) is a $G$-homotopy equivalence, and \(P \to B\) is
a $G$-fibration, then the pullback $G$-fibration $h^*P$ is
$G$-homotopy equivalent to $P$. Applying this
straighforward corollary of Theorem~\ref{theorem:homotopy-pullback}
twice, once on each left-hand square of the two diagrams
in~\eqref{eq:full-pullback}, we may conclude that $Q$ is $G$-homotopy
equivalent to $Q_f$, and $Q'$ is $G$-homotopy equivalent to
$Q'_f$. Since we have already seen that $Q_f$ and $Q'_f$ are
$G$-homotopy eequivalent we have shown that $Q$ and $Q'$ are
$G$-homotopy equivalent, as desired. 
\end{proof}

\bibliographystyle{habbrv}
\bibliography{ref}

\def\cprime{$'$}
\begin{thebibliography}{10}

\bibitem{AMM98}
A.~Alekseev, A.~Malkin, and E.~Meinrenken.
\newblock Lie group valued moment maps.
\newblock {\em J.\ Diff.\ Geom.}, 48:445--495, 1998, dg-ga/9707021.

\bibitem{AtiSeg04}
M.~Atiyah and G.~Segal.
\newblock Twisted {$K$}-theory.
\newblock {\em Ukr. Mat. Visn.}, 1(3):287--330, 2004.

\bibitem{AtiBot83}
M.~F. Atiyah and R.~Bott.
\newblock The {Y}ang-{M}ills equations over {R}iemann surfaces.
\newblock {\em Philos. Trans. Roy. Soc. London Ser. A}, 308(1505):523--615,
  1983.

\bibitem{BTW04}
R.~Bott, S.~Tolman, and J.~Weitsman.
\newblock Surjectivity for {H}amiltonian loop group spaces.
\newblock {\em Invent. Math.}, 155(2):225--251, 2004, math.DG/0210036.

\bibitem{Gol01}
R.~F. Goldin.
\newblock The cohomology rings of weight varieties and polygon spaces.
\newblock {\em Adv. Math.}, 160(2):175--204, 2001, math.SG/0201138.

\bibitem{HarLan07}
M.~Harada and G.~D. Landweber.
\newblock Surjectivity for {H}amiltonian {$G$}-spaces in {$K$}-theory.
\newblock {\em Trans. Amer. Math. Soc.}, 359:6001--6025, 2007, math.SG/0503609.

\bibitem{HL-kernel}
M.~Harada and G.~D. Landweber.
\newblock The ${K}$-theory of abelian symplectic quotients.
\newblock {\em Math. Res. Lett.}, to appear.

\bibitem{JK95}
L.~C. Jeffrey and F.~C. Kirwan.
\newblock Localization for nonabelian group actions.
\newblock {\em Topology}, 34:291--327, 1995.

\bibitem{Kir84}
F.~Kirwan.
\newblock {\em Cohomology of quotients in symplectic and algebraic geometry},
  volume~31 of {\em Mathematical Notes}.
\newblock Princeton University Press, Princeton, N.J., 1984.

\bibitem{May75}
J.~P. May.
\newblock Classifying spaces and fibrations.
\newblock {\em Mem. Amer. Math. Soc.}, 1(1, 155):xiii+98, 1975.

\bibitem{May96}
J.~P. May.
\newblock {\em Equivariant homotopy and cohomology theory}, volume~91 of {\em
  CBMS Regional Conference Series in Mathematics}.
\newblock Published for the Conference Board of the Mathematical Sciences,
  Washington, DC, 1996.

\bibitem{Mil63}
J.~Milnor.
\newblock {\em Morse theory}.
\newblock Based on lecture notes by M. Spivak and R. Wells. Annals of
  Mathematics Studies, No. 51. Princeton University Press, Princeton, N.J.,
  1963.

\bibitem{PS86}
A.~Pressley and G.~Segal.
\newblock {\em Loop groups}.
\newblock Oxford Mathematical Monographs. The Clarendon Press Oxford University
  Press, New York, 1986.

\bibitem{Seg68}
G.~Segal.
\newblock Equivariant {$K$}-theory.
\newblock {\em Inst. Hautes \'Etudes Sci. Publ. Math.}, 34:129--151, 1968.

\bibitem{Sel97}
P.~Selick.
\newblock {\em Introduction to homotopy theory}, volume~9 of {\em Fields
  Institute Monographs}.
\newblock American Mathematical Society, Providence, RI, 1997.

\bibitem{TW03}
S.~Tolman and J.~Weitsman.
\newblock The cohomology rings of symplectic quotients.
\newblock {\em Comm. Anal. Geom.}, 11(4):751--773, 2003, math.DG/9807173.

\bibitem{VeVi03}
G.~Vezzosi and A.~Vistoli.
\newblock Higher algebraic {$K$}-theory for actions of diagonalizable groups.
\newblock {\em Invent. Math.}, 153(1):1--44, 2003, math.AG/0107174.

\end{thebibliography}

\end{document}